\documentclass[12pt]{article}
\usepackage{amsfonts}
\usepackage{amsmath,amssymb}
\usepackage{mathrsfs}
\usepackage{amsxtra}
\usepackage{amstext}
\usepackage{amscd}
\usepackage{latexsym}
\usepackage{bbold}
\newtheorem{thm}{Theorem}[section]
\newtheorem{prop}[thm]{Proposition}
\newtheorem{lem}[thm]{Lemma}
\newtheorem{cor}[thm]{Corollary}
\newtheorem{remark}[thm]{Remark}

\numberwithin{equation}{section}

\def\Tr{{\rm Tr}}

\def\Bbb{\mathbb}

\def\C{\Bbb C}
\def\Z{\Bbb Z}
\def\N{\Bbb N}

\usepackage{color}

\newcommand{\red}{\textcolor{red}}

\newcommand{\thedate}

\begin{document}
\hfill \textbf{\red{DSv35}}

\vspace{-0.5cm}

\[
\mbox{\bf\large {Dynamical Semigroups for Unbounded Repeated}}
\]
\[
\mbox{\bf\large {Perturbation of Open System}}
\]

\begin{center}

\setcounter{footnote}{0}
\renewcommand{\thefootnote}{\arabic{footnote}}

\textbf{Hiroshi Tamura} \footnote{tamurah@staff.kanazawa-u.ac.jp}\\
      {Institute of Science and Engineering}      \\
              and \\
        Graduate School of the Natural Science and Technology\\
           Kanazawa University,\\
          Kanazawa 920-1192, Japan

\medskip

\textbf{Valentin A.Zagrebnov }\footnote{Valentin.Zagrebnov@univ-amu.fr}\\
Institut de Math\'{e}matiques de Marseille - UMR 7373 \\
CMI-AMU, Technop\^{o}le Ch\^{a}teau-Gombert\\
39, rue F. Joliot Curie, 13453 Marseille Cedex 13, France\\
and \\
D\'{e}partement de Math\'{e}matiques \\
Universit\'{e} d'Aix-Marseille - Luminy, Case 901\\
163 av.de Luminy, 13288 Marseille Cedex 09, France

\vspace{0.5cm}

ABSTRACT

\end{center}

\vspace{- 0.5cm}
We consider dynamical semigroups with unbounded Kossakowski-Lindblad-Davies
generators which are related to evolution of an open system with a tuned repeated harmonic perturbation.
Our main result is the proof of existence of uniquely determined minimal trace-preserving
strongly continuous dynamical semigroups on the space of density matrices.
The corresponding dual $W^*$-dynamical system is shown to be unital quasi-free and
completely positive automorphisms of the CCR-algebra.
We also comment on the action of dynamical semigroups on quasi-free states.
\vspace{-0.5cm}

\tableofcontents
\section{Introduction: Repeated Perturbation}\label{Intro}

Since repeated perturbation of Hamiltonian dynamics is piecewise constant, its analysis reduces to study of Quantum
Dynamical Semigroups (QDS) on the space of states and of their generators.
A similar reduction is also valid for
repeated perturbation of open quantum dynamical systems, which are described by
dissipative extensions of Hamiltonian generators \`{a} la Kossakowski-Lindblad-Davies (KLD) \cite{AJP3, BJM} in Markov approximation \cite{AJP2}.
The theory of QDS is quite satisfactory for {\textit{bounded}} generators and for their \textit{bounded} KLD extensions
\cite{Da1}.
A generalisation of this theory to the case of \textit{unbounded} dissipative generators was initiated in \cite{Da2, Da3}
and developed in \cite{Fa, EL, DVV, Pu} for completely positive maps on CCR-algebras.
The progress in construction of the \textit{minimal} dynamical semigroups for
unbounded dissipative generators is essentially due to ideas that come back
to T.Kato \cite{Ka1}.
These ideas were developed first in \cite{Da2}.
Later they inspired the construction  and the abstract analysis of uniqueness and
trace-preserving (or Markovian) property of the minimal QDS with KLD generators,
see \cite{ChF}, \cite{AJP3} Lecture 3.
They were followed by important works in the study of QDS (see, e.g. \cite{FR}),
including some recent analysis of singular (relative bound equals to one)
perturbations of positive and substochastic semigroups on the normal states
\cite{M-K} and abstract  spaces of states. \cite{ALM-K}

This paper is addressed to these problems for unbounded generators by a concrete
quantum dynamics which needs different approach from the works mentioned above.
Our model is a dissipative KLD extension of Hamiltonian dynamical system \cite{TZ}, which gives an open system for boson reservoir.
And its generator is unbounded with the relative bound equals to one.

Our main results are the following:

\noindent We construct the generator of the minimal QDS corresponding to the standard
KLD extension of Hamiltonian dynamical system \cite{TZ}, in Section 2.
We prove in Theorem \ref{generate} that  it generates
\textit{strongly continuous}, \textit{positive}, \textit{contraction} and \textit{trace-preserving} semigroups, i.e. the Markov Dynamical Semigroup (MDS) on the space of trace-class operators.
In Section 3, we establish the explicit formulae for the action of its \textit{dual}
MDS on the Weyl CCR-algebra (Theorem \ref{one-step-dualT}).
This allows to prove that the dual MDS is \textit{completely positive}  (Theorem \ref{MDS-CP}).
Finally we prove that the MDS maps the space of \textit{quasi-free} states
into itself, see Proposition \ref{QF}.

In the rest of this section, we briefly review the model of Hamiltonian dynamics
of \cite{TZ} and recall the standard KLD extension to the open system with a linear boson reservoir.
Then we give a formal definition of the generator for our model for the open
system with repeated interaction.

\medskip

Let  $a$ and $a^*$ be the annihilation and the creation operators defined in the Fock space $\mathscr{F}$ generated by a cyclic vector $\Omega $.
That is, the Hilbert space $\mathscr{F}$ contains the
algebraic span $\mathscr{F}_{\mbox{\tiny fin}}$  of vectors $\{(a^{*})^m \Omega\}_{m\geq 0}$ as a dense subset
and $a, a^*$ satisfy the Canonical Commutation Relations (CCR)
\[
      [a, a^*] = \mathbb{1}, \quad  [a, a] = 0, \quad [a^*, a^*] = 0
      \quad \mbox{on} \quad \mathscr{F}_{\mbox{\tiny fin}}.
\]
We denote by $\{\mathscr{H}_k\}_{k=0}^{N}$ the copies of $\mathscr{F}$ for an arbitrary but finite $N \in \mathbb{N}$ and by
$\mathscr{H}^{(N)} $ the Hilbert space tensor product of these copies:
\begin{equation}\label{H-space}
\mathscr{H}^{(N)}  = \bigotimes_{k=0}^{N} \mathscr{H}_k
\end{equation}
and by $ \Omega_{F}: = \Omega^{\otimes(N+1)} $, its cyclic  vector.

In this space, we define the annihilation and the creation operators
\begin{eqnarray}\label{k-bosons}
b_k :=
\mathbb{1} \otimes  \ldots \otimes \mathbb{1} \otimes a \otimes \mathbb{1} \otimes  \ldots \otimes \mathbb{1} \,  , \ \
b_{k}^* :=
\mathbb{1} \otimes  \ldots \otimes \mathbb{1} \otimes a^* \otimes \mathbb{1} \otimes  \ldots \otimes \mathbb{1}
\end{eqnarray}
for $k = 0, 1, 2, \ldots, N$,  where  the operator $a$, or $a^*$, is the
$(k+1)$-th factor.
On algebraic tensor product $\mathscr{H}_{\mbox{\tiny fin}}^{(N)} :=
\mathscr{F}_{\mbox{\tiny fin}}^{\otimes(N+1)}$, these unbounded operators satisfy the CCR:
\begin{equation}\label{CCR}
[b_k, b^*_{k^\prime}] = \delta_{k,k^\prime} \mathbb{1} , \quad
[b_k, b_{k^\prime}] =  [b^*_k, b^*_{k^\prime}] = 0 \  \
 (k,k^\prime =0, 1, 2, \ldots, N) \ .
\end{equation}

We consider the Hamiltonian of the system with time-dependent \textit{repeated} \textit{harmonic} perturbation \cite{TZ}:
\begin{eqnarray}\label{Ham-Model}
H_{N}(t) = E \, b_0^*b_0 + \epsilon \sum_{k=1}^{N} b_k^*b_k +
\eta \, \sum_{k=1}^{N}\chi_{[(k-1)\tau, k\tau)}(t)\,  (b_0^*b_k +  \ b_k^*b_0) \ ,
\end{eqnarray}
for $t\in [0, N\tau)$, where  $\tau, E,\epsilon , \eta >0 $ and $\chi_{[x,y)}(\cdot)$
is the characteristic function  of semi-open intervals $[x,y)\subset \mathbb{R}$.
Here (\ref{Ham-Model}) denotes the self-adjoint operator on the dense domain
\begin{equation}\label{domH}
           \mathcal{D}_0 := \bigcap_{k=0}^{N}{\rm{dom}}(b^*_{k}b_{k}) \
              \subset \mathscr{H}^{(N)}        \, .
\end{equation}

The model (\ref{Ham-Model}) describes the system $\mathcal{S} + \mathcal{C}_N $,
where subsystem $\mathcal{S}$ corresponding to the kinetic term
$E \, b_0^*b_0$ of the Hamiltonian is repeatedly interacting with a long time-equidistant \textit{chain}
$\mathcal{C}_N = \mathcal{S}_1 + \mathcal{S}_2 + \cdots + \mathcal{S}_N $
of subsystems corresponding to the kinetic terms
$\epsilon \sum_{k=1}^N\, b_k^*b_k$.
The Hilbert space $\mathscr{H}_0$ corresponds to the subsystem $\mathcal{S}$ and
the Huberto space $\mathscr{H}_k$ to the subsystem $\mathcal{S}_k \; $ ($k=1, \cdots, N$).
This visualisation is motivated by a number of physical models, see \cite{BJM}, \cite{NVZ}.

For $t\in[(n-1)\tau, n\tau)$, only subsystem $\mathcal{S}_{n}$ interacts with
$\mathcal{S}$ and the system
$\mathcal{S}+\mathcal{C}_N$ is \textit{autonomous} on this time-interval with the self-adjoint Hamiltonian
\begin{equation}\label{Ham-n}
H_n  =  E \, b_0^*b_0 + \epsilon\sum_{k=1}^{N}b_k^*b_k  + \eta \, (b_0^*b_n +  b_n^*b_0)
\end{equation}
on domain $\mathcal{D}_0$. To keep the operator (\ref{Ham-n}) lower semi-bounded,
we assume that parameters $ E, \epsilon, \eta $ satisfy the condition
\begin{equation}\label{H4}
\textbf{(H1)} \hspace{6cm} \eta^2 \leqslant E \, \epsilon \ . \hspace{6cm}
\end{equation}

We denote by $\mathfrak{C}_{1}(\mathscr{H}^{(N)})$ the Banach space of trace-class operators on $\mathscr{H}^{(N)}$ with trace norm $\|\cdot\|_1$.
Its dual space is isometrically isomorphic to the space of bounded operators on $\mathscr{H}^{(N)}$: $\mathfrak{C}_{1}^{\ast}(\mathscr{H}^{(N)})\simeq \mathcal{L}(\mathscr{H}^{(N)})$.
We consider the dual pair corresponding to the bilinear functional
\begin{equation}\label{dual-Tr}
\langle \phi \, |  A \rangle_{\mathscr{H}^{(N)}} = \Tr_{\mathscr{H}^{(N)}}
        (\phi \, A)  \, ,
\quad \mbox{for} \ (\phi, A) \in
\mathfrak{C}_{1}(\mathscr{H}^{(N)}   ) \times \mathcal{L}(\mathscr{H}^{(N)}) \ .
\end{equation}
Positive operators in $\mathfrak{C}_{1}(\mathscr{H}^{(N)})$ with \textit{unit} trace are called
\textit{density matrices}.
For each density matrix $\rho$, we consider the normal state $\omega_{\rho}(\cdot)$ on $\mathcal{L}(\mathscr{H}^{(N)})$ defined by
\begin{equation}\label{state-S}
\omega_{\rho}(\, \cdot \,)  = \langle \rho \, | \, \cdot \, \rangle_{\mathscr{H}^{(N)}}  \, .
\end{equation}

To describe evolution of the open system corresponding to (\ref{Ham-Model}),
we consider the Kossakowski-Lindblad-Davies (KLD) dissipative extension of the Hamiltonian
dynamics to non-Hamiltonian \textit{master equation}: $\partial_{t}\rho(t) = L_{\sigma}(t)(\rho(t))$, with the time-dependent
generator
\begin{equation}\label{K-L-D-Generator}
L_{\sigma}(t)(\rho) := -i \, [H_N(t),\rho]   + \, \mathcal{Q}(\rho) -
\frac{1}{2} ({\mathcal{Q}}^{\ast}(\mathbb{1}) \rho + \rho \, {\mathcal{Q}}^{\ast}
(\mathbb{1})) \, ,
\end{equation}
for $ t \in [0, N\tau)  $.  \cite{AJP3, AF}
The operator $\mathcal{Q}$ acts on $\rho$ as
\begin{equation}\label{K-L-D1}
      \mathcal{Q}(\rho) = \sigma_{-} \, b_{0} \, \rho \, b^* _{0}
               + \sigma_{+} \, b^*_{0}\, \rho \, b_{0} \ .
\end{equation}
Its dual operator ${\mathcal{Q}}^{\ast}$ is defined by the relation
$\langle \mathcal{Q}(\rho) \, |  A \rangle_{\mathscr{H}^{(N)}}
  = \langle \rho \, | {\mathcal{Q}}^{\ast}(A) \rangle_{\mathscr{H}^{(N)}}$:
\begin{equation}\label{K-L-D2}
{\mathcal{Q}}^{\ast}(A) = \sigma_{-} \, b^{*}_{0} \, A \, b _{0}
        + \sigma_{+} \, b_{0}\, A \, b^*_{0} \, .
\end{equation}

Since the Hamiltonian part of the dynamics is \textit{piecewise} autonomous, the generator
(\ref{K-L-D-Generator}) for $t\in[(k-1)\tau,k\tau)$, $k= 1,2, \ldots , N$, gets the form
\begin{align}\label{Generator-KLD}
L_{\sigma,k}(\rho):= -i[H_k,\rho] +  \, \mathcal{Q}(\rho) -
\frac{1}{2} ({\mathcal{Q}}^{\ast}(\mathbb{1}) \rho + \rho {\mathcal{Q}}^{\ast}(\mathbb{1})) \  .
\end{align}

Note that the form of generators (\ref{K-L-D-Generator}), (\ref{Generator-KLD}) corresponds to repeated perturbation of
the open system $\mathcal{S} + \mathcal{R}$, i.e. we study $(\mathcal{S} + \mathcal{R}) + \mathcal{C}_N$ for external
boson reservoir $\mathcal{R}$. Then a formal solution $\rho(t)$ of the Cauchy problem for the master equation corresponding to
initial condition $\rho(0) = \rho$, is defined by the evolution map $\{T_{t,0}^{\sigma}\}_{t\geq 0}$. It is a composition of
QDS with generators (\ref{Generator-KLD}):
\begin{equation}\label{Cauchy}
\rho(t) = T_{t,0}^{\sigma} (\rho) := (T_{n, \nu(t)}^{\sigma} \, T_{n-1}^{\sigma} \ldots T_{2}^{\sigma} \
T_{1}^{\sigma}) (\rho)
\end{equation}
for $t = (n-1)\tau + \nu(t)$ and $n\leqslant N$, where
 $T_{k,s}^{\sigma} = e^{s L_{\sigma,k}} $, $ T_{k}^{\sigma} = T_{k,\tau}^{\sigma}$
$\   (k = 1, 2, \cdots, n)$.
Consequently, the analysis of evolution for repeated perturbation reduces
to the study of QDS on the intervals $[(k-1)\tau,k\tau), \ k =1, \cdots, N$.

It is known that for the \textit{standard} KLD generator of the form
(\ref{Generator-KLD}) with bounded $H_k$, $\mathcal{Q}$ and
${\mathcal{Q}}^{\ast}$, the corresponding QDS  $\{T_{k,s}^{\sigma}\}_{s\geq 0}$ on
$\mathfrak{C}_{1}(\mathscr{H}^{(N)})$ is norm-continuous, completely positive and trace-preserving, see e.g.\cite{Da1}.
The first aim of the present paper is to give a rigorous meaning to the generator
of the standard form (\ref{Generator-KLD}) with {unbounded} operators (\ref{Ham-n}), (\ref{K-L-D1}) and (\ref{K-L-D2}) and to construct QDS for the solution
(\ref{Cauchy}).
And then, we show the above properties for our QDS with unbounded generators.

\medskip

Our next hypothesis demands that the parameters $\sigma_{\pm}$ (\ref{K-L-D1}), (\ref{K-L-D2}) satisfy the condition:
\begin{equation}\label{H-sigma}
\textbf{(H2)} \hspace{6cm} 0 \leqslant \sigma_{+} < \sigma_{-} \ . \hspace{5cm}
\end{equation}
Together with \textbf{(H1)}, the condition plays an important role in the
construction of semigroups $\{T_{k,s}^{\sigma}\}_{s\geq 0}$ with trace-preserving
property. cf. Theorem \ref{generate}.
Under these hypothesis, complete positivity of the dual semigroups $\{
T_{k,s}^{\sigma \, \ast}\}_{s\geq 0} $ are established in Section \ref{CP-MDS}.

Finally, from now on we suppress the superscript $N$ in $\mathscr{H}^{(N)}$ for brevity.

\section{Minimal Dynamical Semigroup}\label{MED}
\subsection{{Unbounded generators}}\label{UnbGen}

First, we define operators related to the Hamiltonian (\ref{Ham-n}) in the Hilbert space $\mathscr{H}$ (\ref{H-space}):
\begin{eqnarray}
&& K_{0} = \frac{\sigma_+}{2} b_0 b_0^* + \frac{\sigma_-}{2} b_0^* b_0
+i \, \big( (E - \epsilon) b_0^* b_0 + \epsilon \, \hat{n} \big) \, , \  \ \hat{n} = \sum_{j=0}^N b_j^*b_j \, ,
 \label{K0} \\
&& K_{n} = K_0 + i \, \eta (b_0^* b_{n} + b_{n}^* b_0 )
      =  \frac{1}{2} \, {\mathcal{Q}}^{\ast}(\mathbb{1}) + i\, H_{n} \, , \  n =1,2, \ldots, N \ .  \label{Kn}
\end{eqnarray}
Here $E, \epsilon, \eta >0$ and $\sigma_{\pm}$ satisfy  \textbf{(H1)} and \textbf{(H2)}, respectively.
Domains of these operators are identical to $\mathcal{D}_0 $ (\ref{domH}), which is dense in $\mathscr{H}$.

\begin{lem}\label{m-accr} For $n =1,2,\ldots, N$, the operator $K_{n}$ is $m$-accretive.
\end{lem}
\vskip-2mm

\noindent For the proof, see Appendix

\medskip

It is known that for any $m$-\textit{accretive} $A$ in a Hilbert space, the operator
$(-A)$ is the generator of a one-parameter \textit{Strongly Continuous Contraction
Semigroup} (SCCS) $\{e^{-tA}\}_{t\geqslant 0} \, $ on the Hilbert space, in general,
e.g. \cite{Ka2}, \cite{Za}.
Then Lemma \ref{m-accr} implies:
\begin{cor}\label{m-accr-cor}
The operator $-K_n$ is the generator of a SCCS $\{e^{- t \, K_n}\}_{t\geqslant 0}$ on $\mathscr{H}$ for $ n=1, 2, \cdots, N$.
\end{cor}

Next we make precise definition of operators (\ref{Generator-KLD}). Since the operators $\{b_n, b_n ^*\}_{n=0}^{N}$
in $\mathscr{H}$ are unbounded, the operators (\ref{Generator-KLD}) in the Banach space $\mathfrak{C}_1(\mathscr{H})$ are
also unbounded.
Let $\Phi: \mathfrak{C}_1(\mathscr{H}) \rightarrow \mathfrak{C}_1(\mathscr{H})$ be the positive injection defined by
$\Phi(\rho) = (\mathbb{1}+\hat n)^{-1} \rho(\mathbb{1} + \hat n)^{-1}$, and put
$\widetilde{\mathscr{D}} = \Phi ( \mathfrak{C}_1(\mathscr{H})) $.
Note that $\hat n$ is a non-negative self-adjoint operator on domain $\mathcal{D}_0$.
In fact,
\begin{equation}\label{CONS}
\psi_{m} = \frac{b_0^{* m_0} b_1^{* m_1}  \ldots  b_N^{* m_N}} {\sqrt{m_0! \, m_1!\, \ldots \, m_N!}} \
\Omega_{F}
\end{equation}
is the eigenvector of $\hat n$ with eigenvalue $ \sum_{k=0}^N m_k $ for $ m = (m_0, \cdots, m_N) \in \Z_+^{N+1}$.
And the set of vectors (\ref{CONS}) for $m \in \Z_+$ form  a \textit{Complete Ortho-Normal System} (CONS) of $\mathscr{H}$.

Note that operators (\ref{Kn}) are relatively bounded with respect to $(\mathbb{1}+\hat n)$, i.e.,
$\|K_{n} \psi\| \leq \alpha \ \|(\mathbb{1}+\hat n)\psi\|$, $\psi \in \mathcal{D}_0$ hold for some
$\alpha >0$,  \cite{Ka2}.
Taking into account that operators $b_0(\mathbb{1} +\hat n)^{-1}$ and $b_0^*(\mathbb{1} +\hat n)^{-1}$ are bounded,
the unbounded operator (\ref{Generator-KLD})
\begin{equation}\label{defL-0}
L_{\sigma, n}(\rho) = -K_n  \rho -\rho K_n  ^*
               + \sigma_- b_0\rho b_0^* + \sigma_+ b_0^*\rho b_0
\end{equation}
can be defined on $\widetilde{\mathscr{D}}$ as
\begin{equation*}
L_{\sigma,n}(\Phi(\rho)) = -K_{n}(\mathbb{1} +\hat n)^{-1}\rho(\mathbb{1} +\hat n)^{-1}
- (\mathbb{1} +\hat n)^{-1}\rho (K_{n} (\mathbb{1} +\hat n)^{-1})^*
\end{equation*}
\begin{equation}\label{defL}
 + \sigma_-b_0(\mathbb{1} +\hat n)^{-1}\rho (b_0(\mathbb{1} +\hat n)^{-1})^*
 + \sigma_+b_0^*(\mathbb{1} +\hat n)^{-1}\rho (b_0^*(\mathbb{1} +\hat n)^{-1})^* \ ,
\end{equation}
for any $\rho \in \mathfrak{C}_1(\mathscr{H})$ and $n= 1,2 \ldots, N$.
Note that the domain $\widetilde{\mathscr{D}}$ is dense in $\mathfrak{C}_1(\mathscr{H})$, since it contains all finite-rank operators
made of vectors lie in $\mathcal{D}_0$.

\subsection{{Dynamical semigroup on the space of density matrices}}\label{DS-on-DensMatr}
To construct dynamical semigroups (DS) with the generators which are extensions of (\ref{defL-0}),
we recall some results of the Kato-Davies approach \cite{Ka1}, \cite{Da2}.
Since these results are applicable to any $ n = 1, 2, \ldots N$ verbatim,
we describe them under our notations for the case $n=1$ and the corresponding semigroup.

First, we note that the operator $K_1$ (\ref{Kn}) satisfies the identity
\begin{equation}\label{Ident-K-Q}
                 - (K_1 \varphi, \psi)      - (\varphi, K_1  \psi)
       + \sigma_-(b_0\varphi, b_0\psi)   + \sigma_+(b_0^*\varphi, b_0^*\psi) = 0 \ ,
\end{equation}
for all $\varphi, \psi \in \mathcal{D}_0$ (\ref{domH}).

Let $V$ denote the Banach subspace of all \textit{self-adjoint} elements of $\mathfrak{C}_1(\mathscr{H})$.
The family of maps
\begin{equation} \label{St}
S_t(\rho) = e^{-t K_1}\, \rho\, (e^{-t K_1})^*  \   \ (t\geqslant 0 \ , \  \rho \in V)
\end{equation}
defines a positive SCCS on $V$.
Let $Z$ be the generator of $S_t $ and dom$\, (Z)$ its domain.
Then $\mathscr{D} = \Psi(V):=(\mathbb{1}+K_1)^{-1}V((\mathbb{1}+K_1)^{-1})^{\ast} \subset $ dom$\, (Z)$ and
\begin{equation}\label{Z-rho}
Z(\rho) = -K_1\rho - \rho K_1^* \quad  \mbox{ for } \quad \rho \in \mathscr{D}
\end{equation}
hold.
The set $\mathscr{D}$ is dense in $V$ and  a \textit{core} of the generator $Z$.
Note that $\mathscr{D} = \widetilde{\mathscr{D}}\cap V$.
There are two positive $Z$-bounded operators $J_-$ and $J_+$ on dom$\, Z$
such that
\begin{equation}\label{Jpm}
J_-(\rho) = b_0\rho b_0^* \ ,  \quad  J_+(\rho) = b_0^*\rho b_0 \qquad \mbox{ for } \quad \rho \in \mathscr{D} \,.
\end{equation}
Then, the operator $\hat{L}:= Z +  \sigma_- J_-  + \sigma_+ J_+ $ is defined
on the domain dom$\, (Z)$.
Whereas let us denote $L$ the operator
(\ref{defL-0}) for $n=1$ with domain  $\widetilde{\mathscr{D}}$.
Here we understand (\ref{Z-rho}) and  (\ref{Jpm}) as in (\ref{defL}).
Then,
\begin{equation}\label{Ident}
\Tr_{\mathscr{H}}(\hat{L}(\rho)) = 0  \quad \mbox{ holds for } \quad
\rho \in {\rm dom}\, (Z)
\end{equation}
and the operator $J:=(\sigma_- J_- + \sigma_+J_+)$ is $Z$-bounded with the {relative} bound equals to \textit{one},
which require non-perturbative arguments to construct the DS corresponding to $\hat{L}$.
\begin{prop}\label{Tr}
For any $r \in [0, 1)$ the operator $Z+ r (\sigma_- J_- + \sigma_+ J_+ )$ with domain
${\rm dom}\,(Z)$ is the generator of a positive SCCS $\{T_{t,r}\}_{t \geqslant 0}$ on $V$.
\end{prop}

\begin{prop}\label{minT}
There exists a positive SCCS $\{T_t\}_{t \geqslant 0}$ on $V$ such that
\[
                  \lim_{r \to 1} T_{t,r}(\rho) =T_t(\rho) \ , \ \rho \in V \ ,
\]
uniformly in each compact interval of $t\geqslant 0$. The generator $M$ of $T_t$ is a closed extension of
the operator $\hat{L}$.
\end{prop}
\begin{remark}\label{minT-b} Since perturbation $J$ has relative bound $1$, the operator $\hat{L}$ may have many closed
extensions \cite{Ka2}. The semigroup constructed in Proposition \ref{minT} is
minimal in the following sense: if the SCCS $\{T_{t}^{\prime}\}_{t \geqslant 0}$ has the generator $M^{\prime}$, which is
another extension of $\hat{L}$, then $ T_{t}^{\prime} > T_{t} $ holds for all $ t > 0$.
Moreover, in spite of (\ref{Ident-K-Q}), or the "conservativity" (\ref{Ident}),
the minimal DS need not be trace-preserving.
\end{remark}
\begin{prop}\label{coreM}
If ${\rm dom}\,(Z)$ is a core of the generator $M$, then the minimal semigroup
$\{T_t\}_{t \geqslant 0}$ is trace-preserving, i.e. a Markovian semigroup.
\end{prop}
\medskip

Thus far, we have got a glimpse of results from \cite{Ka1, Da2}.
Now we come back to analysis of our concrete open system (\ref{Ham-n}), (\ref{Generator-KLD}) for the master equation with
generators $L_{\sigma,n}$ (\ref{defL-0}) on domain $\widetilde{\mathscr{D}}$.
\begin{thm}\label{generate}
For each $n = 1,2 \ldots, N$, the closure of the operator $L_{\sigma,n}\!\upharpoonright\!_{\mathscr{D}}$
is the generator of a trace-preserving SCCS on $V$.
\end{thm}
\noindent {\sl Proof:}
It is enough to consider only the case $n =1$ as above.

\noindent$1^{\circ \ }$ We start by checking that ${\rm dom}\,(Z)$ is a core of $M$. Let us define
the SCCS $R_s$ on $V$ by
\[
  R_s(\rho) = e^{-s \hat n}\rho e^{-s \hat n} \ , \ \  ( s \geqslant 0 \ , \ \rho \in V )\, .
\]
Then
\begin{equation} \label{CR11}
    e^{-tK_1}e^{-s\hat n} = e^{-s\hat n} e^{-tK_1} \qquad \mbox{on } \; \mathscr{H} ,
\end{equation}
as well as $b_0e^{-s\hat n} = e^{-s}e^{-s\hat n}b_0 $, $b_0^*e^{-s\hat n} = e^{s}e^{-s\hat n}b_0^* $ on $\mathcal{D}_0$.
Combining with (\ref{St}), we obtain
\begin{equation}\label{CR21}
R_s(S_t(\rho)) = S_t(R_s(\rho)) \qquad \mbox{for } \; \rho \in V \, .
\end{equation}

Since any element of $\mathscr{D}$ can be expressed as a convergent sum of the rank-one operators with eigenvectors in
$\mathcal{D}_0$, we obtain that
\begin{eqnarray}  \notag
        J_+(R_s(\rho)) &=& e^{2s}R_s(J_+(\rho)) \,  ,
\\  \label{CR22}
         J_-(R_s(\rho)) &=& e^{-2s}R_s(J_-(\rho)) \, ,
\end{eqnarray}
hold for $\rho \in \mathscr{D}$. Differentiating (\ref{CR21}) with respect to $t$, one gets
\begin{equation}\label{R1}
R_s({\rm dom}\,(Z)) \subset {\rm dom}\,(Z) \quad \mbox{and} \quad
R_s(Z(\rho)) = Z(R_s(\rho))  \quad \mbox{for} \quad \rho \in {\rm dom}\,(Z).
\end{equation}
Note that $Z$-boundedness of $J_{\pm}$ together with (\ref{CR22}) and boundedness of
$J_{\pm}R_s $ imply
that the same relations (\ref{CR22}) hold for all $\rho \in {\rm dom}\,(Z)$.
Hence, (\ref{R1}) and (\ref{CR22}) yield
\begin{equation}\label{RZ}
         (Z + \sigma_+ J_+ + \sigma_- J_- )R_s = R_s( Z + e^{2s}\sigma_+ J_+
     + e^{-2s} \sigma_- J_- )
\end{equation}
on dom$\,(Z)$.
Now we  introduce the operators $\tilde K_0$ and $\tilde K_1$ which are defined
by replacing parameters $\sigma_{\pm}, E, \epsilon$ and $\eta$ in $K_0$ and $K_1$
(see (\ref{K0}), (\ref{Kn})) by
$\, \tilde \sigma_{\pm} = e^{\pm 2s}\sigma_{\pm}, \ \tilde E = r(s)E, \ \tilde\epsilon = r(s)\epsilon \, $ and
$\tilde \eta =r(s)\eta$,  where
\begin{equation}\label{r(s)}
                  r(s) := \frac{e^{2s}\sigma_+ + e^{-2s}\sigma_-}{\sigma_+ + \sigma_-}  .
\end{equation}
To keep $r(s) \in (0,1)$, we set $s \in (0, 2^{-1}\log\sigma_-/\sigma_+)$ that is
possible by the hypothesis \textbf{(H2)}: $0 \leqslant \sigma_+ < \sigma_-  $.
Note that $\lim_{s \downarrow 0}r(s) =1$.

By virtue of (\ref{Kn}), one gets the identity
\[
          K_1 = \frac{\tilde K_1}{r(s)} - \frac{\sigma_+\sigma_-}{\sigma_+ + \sigma_-}
               \frac{\sinh 2s}{r(s)} \ \mathbb{1} \ .
\]
Then we obtain that
\begin{equation}\label{ZtildeZ}
       Z = \frac{\tilde Z}{r(s)} + \frac{2\sigma_+\sigma_-}{\sigma_+ + \sigma_-}
        \frac{\sinh 2s}{r(s)} \ \mathbb{1}
\end{equation}
holds on $\mathscr{D}$. Here operator $\tilde Z$ is given by the same expression as (\ref{Z-rho}), but with $\tilde K_1$
instead of $K_1$. Taking the closure in equality (\ref{ZtildeZ}), one gets that
${\rm dom}\, (Z) = {\rm dom}\, (\tilde Z)$ and  that (\ref{ZtildeZ}) holds also on ${\rm dom}\, (Z)$.
Hence, the operators $\tilde{J_{\pm}}$ which are $\tilde Z$-bounded
extension of (\ref{Jpm}), are equal to $J_{\pm}$, respectively.
Therefore, the equality
\[
       Z + \tilde\sigma_+ J_+ + \tilde\sigma_- J_-  = \frac{1}{r(s)}
   \Big[\frac{2\sigma_+\sigma_-}{\sigma_+ + \sigma_-}\sinh 2s \ \mathbb{1} + \tilde Z +
            r(s) (\tilde\sigma_+ J_+ + \tilde\sigma_- J_- ) \Big] \ ,
\]
also holds on ${\rm dom}\, (Z)$.
Together with ({\ref{RZ}), this yields the relation
\begin{eqnarray}\label{R3}
&&( \lambda\mathbb{1} - Z - \sigma_+ J_+ - \sigma_- J_- )\, R_s = \\
&&\frac{1}{r(s)} R_s \Big [\Big(r(s) \lambda - \frac{2\sigma_+\sigma_-}{\sigma_+ + \sigma_-}
         \sinh 2s \Big) \mathbb{1}
- \tilde Z -  r(s) (\tilde\sigma_+ J_+ + \tilde\sigma_- J_- ) \Big] \ .  \nonumber
\end{eqnarray}
on ${\rm dom}\, (Z)$.
Now, for arbitrary $\lambda >0$, we choose $s \in (0, 1)$  small enough such that :
\begin{equation*}\label{range}
r(s) \lambda - \frac{2\sigma_+\sigma_-}{\sigma_+ + \sigma_-}\sinh 2s > 0 \ .
\end{equation*}
Proposition \ref{Tr} in the tilded context yields that
$\tilde Z + r(s) (\tilde\sigma_+ J_+ + \tilde\sigma_- J_- )$ is the generator of a SCCS.
Hence by the Hille-Yosida theorem, its resolvent set includes
$\C_+: = \{ \, z \in \C \, | \, {\rm Re} \, z > 0 \, \}$,
which yields that the last factor in the right-hand side of (\ref{R3}) is invertible
and that the range of the operator in the
left-hand side: $(\lambda\mathbb{1} - Z - \sigma_+ J_+ - \sigma_- J_- )\, R_s \, (V)$
coincides with the set $R_s \, (V)$, which is obviously dense in $V$.
Hence, the range of $(\lambda\mathbb{1} - Z - \sigma_+ J_+ - \sigma_- J_-)$ is also dense
in $V$.

Note that by Proposition \ref{minT} the operator $\hat{L} = Z + \sigma_+ J_+ + \sigma_- J_- $
on the domain ${\rm dom}\,(Z)$
is closable since it has the closed extension $M$. Let $M_0$ be the closure of $\hat{L}$. Then we have:
$\lambda\mathbb{1} - M \supseteq \lambda - M_0 \supset  \lambda\mathbb{1} - \hat{L}$, which implies for $\lambda >0$:
\begin{equation}\label{R4}
(\lambda\mathbb{1} - M)^{-1} \supset  (\lambda\mathbb{1} - M_0)^{-1} \supset (\lambda\mathbb{1} - \hat{L})^{-1} \ .
\end{equation}}
By the conclusion in the previous paragraph, the domain of the last operator in (\ref{R4}) is dense.
Hence, by Proposition \ref{minT} the first operator in (\ref{R4}) is a closed bounded extension of the last.
Since the second operator is the closure of the last one and a restriction
of the bounded $(\lambda\mathbb{1} - M)^{-1}$, then it is also a bounded operator on $V$.
This yields $M = M_0$, which implies that the minimal semigroup is trace-preserving
by Proposition \ref{coreM}.

\medskip

\noindent$2^{\circ \ }$ To finish the proof, we show that $\mathscr{D}$ is a core of $M$. We have already established that
$\mathscr{D}$ is a core of $Z$ and that ${\rm dom}\, (Z)$ is a core of $M_0 = M$. Therefore, for any
$\rho \in {\rm dom}\, (M)$, there exists a sequence $\{ \rho_m \}_{m \geqslant 1} \subset {\rm dom}\, (Z)$ such that
$\rho_m \to \rho, M\rho_m \to M\rho$, as $m \rightarrow \infty$.
Since $ \mathscr{D} $ is a core of $Z$, for each $m$
there exists a sequence $\{ \rho_{m,k} \}_{k \geqslant 1} \subset \mathscr{D}$ such that
$\rho_{m,k} \to \rho_m, Z\rho_{m,k} \to Z\rho_m $  for  $k \rightarrow \infty$. Then by the $Z$-boundedness of $J_{\pm}$, we
also have $J_{\pm} \rho_{m,k} \to J_{\pm} \rho_m $, and thereby $M \rho_{m,k} \to M \rho_m $ .

Therefore, we can choose a diagonal sequence $\{ \rho_{m,k_m} \}_{m\geqslant 1} \subset \mathscr{D}$
such that
\begin{equation}
       \rho_{m,k_m} \to  \rho \quad \mbox{and}   \quad    M \rho_{m,k_m} \to M \rho
            \quad \mbox{hold for $n\rightarrow\infty$}\, .
\end{equation}
Hence, the closure of the operator $L_{\sigma, 1}\!\upharpoonright\!_{\mathscr{D}}$ coincides with $M$.
This completes the proof of the theorem for $n =1$.
\hfill {$\square$}
\begin{remark}\label{rmkex}
The set of density matrices
$\{ \, \rho \in \mathfrak{C}_1(\mathscr{H}) \, | \, \rho \geqslant 0, \Tr_{\mathscr{H}} \, \rho =1 \, \} \subset V$
is obviously invariant subset of $\mathfrak{C}_1(\mathscr{H})$ for the Markov Dynamical Semigroups (MDS)
$\{T_{n,t}^{\sigma}\}_{t \geqslant 0}$, $n = 1,2, \ldots, N$.
On the other hand, the semigroups $ \{T_{n,t}^{\sigma}\}_{t \geqslant 0} $ can be
extended to the MDS on the Banach space $\mathfrak{C}_1(\mathscr{H})$ by linearity.
\end{remark}
\vspace{-0.5cm}
\section{{Markov Dynamical Semigroup on Dual Space}}\label{DDS}
\subsection{{Dual dynamics}}\label{DD}
Equivalent and often more convenient description of the evolution
$\rho \mapsto T_{t,0}^{\sigma}(\rho), \rho \in \mathfrak{C}_1(\mathscr{H})$ is
the dual evolution $\{T_{t,0}^{\sigma \, \ast}\}_{t\geqslant 0}$
on the dual space $\mathfrak{C}_{1}^{\ast}(\mathscr{H})\simeq \mathcal{L}(\mathscr{H})$.

For repeated perturbation, we have to study semigroups $\{T_{n,t}^{\sigma\; *}\}_{t \geqslant 0}$
dual to the SCCS $\{T_{n,t}^{\sigma}\}_{t \geqslant 0}$ constructed
in Theorem \ref{generate}:
\begin{equation}\label{T-*sigma}
\langle T^{\sigma}_{n, t}(\rho)\;|\; A \rangle_{\mathscr{H}}=
\langle \rho \;|\; T^{\sigma \; *}_{n, t}(A)\rangle_{\mathscr{H}} \quad \mbox{for} \quad
(\rho, A) \in \mathfrak{C}_{1}(\mathscr{H}) \times \mathcal{L}(\mathscr{H}) \,, \
n=1, \cdots, N \,.
\end{equation}
Since the maps $T_{n,t}^{\sigma}$ are trace-preserving, the dual semigroups are \textit{unital} (unity-preserving)
contractions. They are also called the Markov Dynamical Semigroups (MDS).

Because the semigroup $\{T_{n,t}^{\sigma}\}_{t \geqslant 0}$ has unbounded generator,
the adjoint semigroup $\{T_{n,t}^{\sigma\, *}\}_{t \geqslant 0}$ is not strongly
continuous on the dual space $\mathcal{L}(\mathscr{H})$.
The duality relation (\ref{T-*sigma}) and the strong continuity of semigroup $\{T_{n,t}^{\sigma}\}_{t \geqslant 0}$
merely imply the weak$^*$-continuity of $T_{n,t}^{\sigma\; *}$ on  $\mathcal{L}(\mathscr{H})$.
Therefore, the pair $( \mathcal{L}(\mathscr{H})  , T_{n,t}^{\sigma\; *})$ is a $W^*$-dynamical system.

Let $\mathscr{A}(\mathscr{H})$ denote the Weyl CCR-algebra on $\mathscr{H}$.
This unital algebra is generated as operator-norm closure of the linear span $\mathscr{A}_{\rm fin}(\mathscr{H})$
of the Weyl operators
\begin{equation}\label{Wz-bis}
W(\zeta) = \exp[i {\big(\langle \zeta, b\rangle + \langle b, \zeta\rangle\big)}/\sqrt 2] \ ,
\end{equation}
where the sesquilinear form notations
\begin{equation}\label{b-forms}
      \langle \zeta, b\rangle  := \sum_{j=0}^N \bar{\zeta}_j b_j,  \qquad
      \langle b, \zeta\rangle  := \sum_{j=0}^N  \zeta_j b^*_j
\end{equation}
are used.
We comment that CCR (\ref{CCR}) has the Weyl form:
\begin{equation}\label{W-CCR}
 W(\zeta_1) W(\zeta_2) = e^{- {i} \, {\rm{Im}}\langle\zeta_{1},
 \zeta_{2}\rangle/2} \ W(\zeta_1 + \zeta_2) \
 \qquad \mbox{for} \qquad \zeta_1,  \zeta_2 \in \C^{N+1} \, .
\end{equation}
and the algebra $\mathscr{A} (\mathscr{H}) $ is dense subset of
$\mathcal{L}(\mathscr{H})$ in the weak as well as in the strong
operator topologies.
(see e.g.  \cite{AJP1} Lectures 4 and 5).

In the rest of this section, we give the explicit form for the action of
$\{{T}_{n,t}^{\sigma \; *}\}_{t\geqslant 0}$ for $1 \leqslant n \leqslant N$
on the Weyl operators.
To this aim, we introduce $(N+1)\times (N+1)$ Hermitian matrices $J_{n}$,
$X_{n}$ and $Y_{n}$ by
\begin{equation}\label{J-n}
      (J_{n})_{jk} = \begin{cases}
                      1 & \quad ( j=k=0 \; \mbox{ or } \; j=k={n}) \\
                      0 & \quad \mbox{otherwise}
                   \end{cases} ,
\end{equation}
\begin{equation}\label{X-n}
     (X_{n})_{jk} = \begin{cases}
                      (E-\epsilon)/2 & \quad (j,k)=(0,0)  \\
                      -(E-\epsilon)/2 & \quad (j,k)=(n, n)  \\
                      \eta & \quad (j,k)=(0,{n}) \\
                      \eta & \quad (j,k)=({n},0) \\
                      0 & \quad \mbox{otherwise}
                   \end{cases}
\end{equation}
and
\begin{equation}\label{Y-n}
    Y_{n}  =\epsilon I + \frac{E-\epsilon}{2}J_{n} + X_{n} \  \qquad
   \mbox{for} \qquad {n}= 1, \cdots, N ,
\end{equation}
where $I$ is the $(N+1)\times (N+1)$ identity matrix. By $P_0$ we denote the $(N+1)\times (N+1)$ matrix:
$(P_0)_{j k} = \delta_{j 0}\delta_{k 0} \ $ $(j, k =1,2, \ldots, N)$. Then the Hamiltonian (\ref{Ham-n}) takes the form
\begin{equation}\label{HY}
   H_{n}={\sum_{j,k=0}^{N} (Y_{n})_{jk}b_j^*b_k} \,.
\end{equation}
\begin{thm}\label{one-step-dualT}
For ${n} = 1, 2, \ldots, N  $ , the action of $\{ {T}_{{n},t}^{\sigma \; \ast}\}_{t\geqslant 0}$
on the Weyl operator has the form:
\begin{equation}\label{T*-n}
{T}_{{n},t}^{\sigma \; \ast}(W(\zeta))  = \Gamma_{{n},t}^{\sigma}(\zeta)
W(U_{n}^{\sigma}(t)\zeta) \ , \  \ \zeta \in \C^{N + 1} \ ,
\end{equation}
\begin{equation}\label{Omega}
\Gamma_{{n},t}^{\sigma}(\zeta) =
\exp \! \Big[ \, - \frac{1}{4} \ \frac{\sigma_- + \sigma_+}{\sigma_- - \sigma_+}
\big(\langle \zeta, \zeta \rangle  - \langle U_{n}^{\sigma}(t)\zeta, U_{n}^{\sigma}(t)\zeta \rangle\big) \Big] \ ,
\end{equation}
and
\begin{equation}\label{U-n-sigma}
U_{n}^{\sigma}(t) = \exp \! \Big[ \, it\Big(Y_{n}  + i \ \frac{\sigma_- - \sigma_+}{2}P_0\Big)\Big] \ .
\end{equation}
\end{thm}
\begin{remark}\label{E-param}
The main effect of non-zero $\sigma_{\mp}$, in comparison to the case
$\sigma_{\mp} = 0$  {\rm{\cite{TZ}}}, may be summarised as an imaginary shift of the energy parameter:
\[
           E \; \to \; E_{\sigma } := E + i \ \frac{\sigma_- - \sigma_+}{2} \ , \ \ \ 0 \leq \sigma_+ < \sigma_- \ .
\]
Note that by {\rm{\textbf{(H2)}}} ${\rm{Im}}(E_{\sigma }) > 0$. Thereby the semigroup $\{U_{n}^{\sigma}(t)\}_{t\geqslant 0}$
is contraction.
\end{remark}
{\sl Proof (of Theorem \ref{one-step-dualT}):} Without loss of generality, we only consider  ${n}=1$.
We put
\begin{equation} \label{Omegazeta}
      \Omega(t) = \Gamma_{1,t}^{\sigma}(\zeta) \quad \mbox{and} \quad \zeta(t) = U_1^{\sigma}(t)\zeta \  .
\end{equation}
$1^{\circ}$ \  The operator-valued equation
\[
       \partial_t(\Omega(t)W(\zeta(t)))   = \Omega(t)\Big(i[H_1,W(\zeta(t))]
     + \sigma_- \ b_{0}^* \ W(\zeta(t))\ b_{0}  \nonumber
\]
\begin{equation}  -\frac{\sigma_-}{2}\{b^*_{0}b_{0},W(\zeta(t))\}
   + \sigma_+ \ b_{0}\ W(\zeta(t)) \ b_{0}^*
     - \frac{\sigma_+}{2} \{b_{0}b^*_{0}, W(\zeta(t))\}\Big) \ .
\label{dual-evol-eq}
\end{equation}
holds on $\mathcal{D}_0$ (\ref{domH}).
Here the derivative in the left-hand side is valid in the strong-operator convergence sense.

To check (\ref{dual-evol-eq}), we use (\ref{HY}) to rewrite the right-hand side as
\[
       \Omega(t)\Big(
      \frac{\langle\zeta(t), Y_1 b\rangle - \langle b, Y_1 \zeta(t)\rangle}
         {\sqrt 2} + i \ \frac{\sigma_+ - \sigma_-}{2}
         \frac{b_0^*\zeta(t)_0 +\overline{\zeta(t)_0}b_0}{\sqrt 2}
\]
\begin{equation}
         - \frac{\sigma_+ +\sigma_-}{4}|\zeta(t)_0|^2
       -i \ \frac{\langle\zeta(t), Y_1 \zeta(t)\rangle}{2}\Big)W(\zeta(t)) \ .
\label{LW}
\end{equation}
Here we took into account the commutators
\[
     [b_k, W(\zeta(t))] = i\ \frac{\zeta(t)_k}{\sqrt 2}W(\zeta(t))\ , \quad
     [b_k^*, W(\zeta(t))] = -i \ \frac{\overline{\zeta(t)_k}}{\sqrt 2} W(\zeta(t)) \ .
\]
On the other hand, using
\[
       \partial_tW(\zeta(t)) =  \Big(i\frac{\langle\partial_t\zeta(t), b\rangle
       + \langle b, \partial_t\zeta(t)\rangle}{\sqrt 2}
\]
\[
          +\frac{1}{2}\Big[i\frac{\langle\zeta(t), b\rangle
       + \langle b, \zeta(t)\rangle}{\sqrt 2} \,  , \,
      i\frac{\langle\partial_t\zeta(t), b\rangle
       + \langle b, \partial_t\zeta(t)\rangle}{\sqrt 2} \Big] \,
       \Big) W(\zeta(t))  \ ,
\]
one gets for the left-hand side of (\ref{dual-evol-eq})
\begin{eqnarray}\label{timeW}
&&\partial_t(\Omega(t)W(\zeta(t))) = \partial_t \Omega(t) W(\zeta(t)) + \\
&& \Big(i\frac{\langle\partial_t\zeta(t), b\rangle + \langle b, \partial_t\zeta(t)\rangle}{\sqrt 2} +
\frac{\langle\partial_t\zeta(t), \zeta(t)\rangle - \langle \zeta(t), \partial_t\zeta(t)\rangle}{4} \Big)\Omega(t) W(\zeta(t))
\ . \nonumber
\end{eqnarray}
Note that expressions (\ref{LW}), (\ref{timeW}), make sense on $\mathcal{D}_0$, which coincides with domain
of the number operator $\hat n \,$. See, e.g., \cite{AJP1}, Lecture 5.

Now the assertion (\ref{dual-evol-eq}) follows, since the equalities
\begin{equation}\label{tzeta}
    \partial_t\zeta(t) = i\Big(Y_1 + i \ \frac{\sigma_- - \sigma_+}{2}P_0\Big)\zeta(t) \, ,
    \quad  \partial_t\Omega(t) = -\frac{\sigma_+ + \sigma_-}{4}|\zeta(t)_0|^2 \, \Omega(t) \,,
\end{equation}
\[
  \frac{\langle \partial_t\zeta(t), \zeta(t)\rangle - \langle \zeta(t), \partial_t\zeta(t)\rangle}{4}
       = - i \frac{\langle\zeta(t), Y_1 \zeta(t)\rangle}{2}
\]
and
\begin{equation*}
   \frac{\langle\partial_t\zeta(t), b\rangle + \langle b, \partial_t\zeta(t)\rangle}{\sqrt 2}
   =  - i \frac{\langle\zeta(t), Y_1 b\rangle - \langle b, Y_1 \zeta(t)\rangle}{\sqrt 2}
  - \frac{\sigma_--\sigma_+}{2}\, \frac{\overline{\zeta(t)_0} \, b_0 + b_0^* \, \zeta(t)}{\sqrt{2}}
\end{equation*}
hold by definitions (\ref{Omegazeta}).

\medskip

\noindent $2^{\circ}$ For any $\rho \in \mathscr{D} = (\mathbb{1} + \hat n)^{-1}V(\mathbb{1} +\hat n)^{-1}$,
the following equality holds:
\begin{equation}\label{diff1}
     \partial_t \, \Tr[\rho \, \Omega(t)W(\zeta(t))] =
                 \Tr[ (L_{\sigma, 1}\rho)\Omega(t)W(\zeta(t))] \, .
\end{equation}

In fact, let $\rho = (\mathbb{1} + \hat n)^{-1} \nu  (\mathbb{1} + \hat n)^{-1}$, where $\nu \in V$ is
approximated  by a family of finite-rank self-adjoint operators $\{ \nu_k \}_{k \geqslant 1}$, i.e.,
$\nu_k \to \nu$, when $k \rightarrow \infty $, in the trace-norm topology.
Then from $1^{\circ}$ and by definition (\ref{defL}), we obtain
\[
 \partial_t\Tr[(\mathbb{1} + \hat n)^{-1} \nu_k  (\mathbb{1} + \hat n)^{-1}\Omega(t)W(\zeta(t))]
\]
\[
        = \Tr \Big[\nu_k\Big\{
  \Big(\Big(-iH_1-\frac{\sigma_-}{2}b_0^*b_0 -\frac{\sigma_+}{2}b_0b_0^*\Big)
      (\mathbb{1} + \hat n)^{-1}\Big)^*\Omega(t)W(\zeta(t))(\mathbb{1} + \hat n)^{-1}
\]
\[
  + (\mathbb{1} + \hat n)^{-1}\Omega(t)W(\zeta(t))
     \Big(\Big(-iH_1-\frac{\sigma_-}{2}b_0^*b_0 -\frac{\sigma_+}{2}b_0b_0^*\Big)
      (\mathbb{1} + \hat n)^{-1}\Big)
\]
\[
       + \sigma_-\big(b_0(\mathbb{1} +\hat n)^{-1}\big)^*
     \Omega(t)W(\zeta(t))\big(b_0(\mathbb{1} +\hat n)^{-1}\big)
\]
\[
    + \sigma_+\big(b_0^*(\mathbb{1}+\hat n)^{-1}\big)^*
          \Omega(t)W(\zeta(t))\big(b_0^*(\mathbb{1} +\hat n)^{-1}\big)\Big\}\Big]
\]
\[
      = \Tr\big[L_{\sigma, 1}\big((\mathbb{1} +\hat n)^{-1}\nu_k(\mathbb{1} +\hat n)^{-1}\big)
       \Omega(t)W(\zeta(t))\big] .
\]
One also gets that the limit:
\[
  \lim_{k\rightarrow\infty} L_{\sigma, 1}\big((\mathbb{1} +\hat n)^{-1} \nu_k (\mathbb{1}+\hat n)^{-1}\big)
     =  L_{\sigma, 1}\big((\mathbb{1} +\hat n)^{-1}\nu(\mathbb{1} +\hat n)^{-1}\big) \, ,
\]
holds in the trace-norm, since by (\ref{defL}) the expression of
$ L_{\sigma, 1}\big((\mathbb{1} +\hat n)^{-1} (\nu_k - \nu) (\mathbb{1} +\hat n)^{-1}\big)$ is the sum of the
products of $\nu_k - \nu$ and $k$-independent bounded operators .
Then
\[
       \partial_t\Tr[(\mathbb{1}+ \hat n)^{-1} \nu_k  (\mathbb{1}+ \hat n)^{-1}\Omega(t)W(\zeta(t))]
         \rightarrow  \Tr\big[L_{\sigma, 1}\big((\mathbb{1} +\hat n)^{-1}\nu(\mathbb{1}+\hat n)^{-1}\big)
       \Omega(t)W(\zeta(t))\big]
\]
holds uniformly in $t$. On the other hand, the limit
\[
        \Tr[(\mathbb{1} + \hat n)^{-1} \nu_k  (\mathbb{1}+ \hat n)^{-1}\Omega(t)W(\zeta(t))]
   \rightarrow  \Tr[(\mathbb{1}+ \hat n)^{-1} \nu  (\mathbb{1}+ \hat n)^{-1}\Omega(t)W(\zeta(t))]
\]
also holds for $k\rightarrow\infty$ uniformly in $t$. Then we obtain the assertion by the standard argument on
differentiation under the limit.

\medskip

\noindent $3^{\circ}$ The equality (\ref{diff1}) also holds for
$\rho \in $ dom$\, \overline{L_{\sigma,1}\!\upharpoonright\!_{\mathscr{D}}}$.
Here $ \overline{L_{\sigma,1}\!\upharpoonright\!_{\mathscr{D}}}$ denotes the closure of the restriction
$ L_{\sigma,1}\!\upharpoonright\!_{\mathscr{D}}$ (c.f. Theorem \ref{generate}).

In fact, for any $\rho \in $ dom$\, \overline{L_{\sigma,1}\!\upharpoonright\!_{\mathscr{D}}}$,
there exists a sequence $\{\rho_k\}_{k\geqslant 1} \subset \mathscr{D}$ such that
\[
     \rho_k \to \rho, \quad  L_{\sigma,1}\!\upharpoonright\!_{\mathscr{D}}
    \, \rho_k \to \overline{L_{\sigma,1}\!\upharpoonright\!_{\mathscr{D}}} \, \rho
\ ,
\]
as $k\rightarrow\infty$, in the trace-norm topology. Then we obtain the assertion  by differentiation
under the limit as in $2^{\circ}$.

\medskip

\noindent $4^{\circ}$ For each $\rho \in $ dom$\, \overline{L_{\sigma,1}\!\upharpoonright\!_{\mathscr{D}}}$ ,
$\zeta \in \C^{N+1}$ and $t\geqslant 0$, the following equality holds:
\begin{equation}\label{TW}
        \Tr\,[T_{1,t}^{\sigma}(\rho)W(\zeta)] =  \Tr\,[ \rho \, \Omega(t)W(\zeta(t))] \ .
\end{equation}

To this aim, we define the function
\[
             f(s, t)  :=  \Tr\,[T_{1,s}^{\sigma}(\rho)\Omega(t)W(\zeta(t))] \quad
                     \mbox{for} \quad s,t \geqslant 0.
\]
Then Theorem \ref{generate} and the Hille-Yosida theorem yield $T_{1,s}^{\sigma}(\rho) \in $
dom$\, \overline{L_{\sigma,1}\!\upharpoonright\!_{\mathscr{D}}}$ and
$\partial_s f( s, t) = \Tr\,[ \overline{L_{\sigma,1}\!\upharpoonright\!_{\mathscr{D}}}
(T_{1,s}^{\sigma}(\rho)) \Omega(t)W(\zeta(t))]$, which is
equal to $\partial_t f( s, t)$ by $3^{\circ}$. Then we obtain $\partial_s f( t - s, s) = 0 $ and
the assertion (\ref{TW}) follows from $f(t,0) = f(0, t)$.

\medskip

\noindent $5^{\circ}$ Since $ T_{1,t}^{\sigma} $ is bounded and dom$\,
\overline{L_{\sigma,1}\!\upharpoonright\!_{\mathscr{D}}}$ is dense in $V$, (\ref{TW}) holds for any $ \rho \in V $.
Note that any $\rho \in \mathfrak{C}_1(\mathscr{H})$ can be presented
as a linear combination of elements from $V$.
The theorem then follows by Remark \ref{rmkex} and by the duality (\ref{T-*sigma}).
\hfill $\square$
\medskip

From (\ref{Cauchy}), the dual evolution map for the repeated perturbation is given by
\begin{equation}\label{T*-t}
T_{t,0}^{\sigma \, *} =T_1^{\sigma \, *} \ \cdots \ T_{n-1}^{\sigma \, *} \
T_{n,\nu(t)}^{\sigma \, *}\, ,
\end{equation}
where $t = (n-1)\tau + \nu(t)$, $n\leqslant N$.
\begin{cor}\label{cor63}
The composition of dual evolutions (\ref{T*-t}) on the Weyl operator is:
\begin{equation*}\label{}
T_{N\tau,0}^{\sigma \; *}(W(\zeta)) = \exp \! \Big[\, - \frac{1}{4}\, \frac{\sigma_- + \sigma_+}{\sigma_- - \sigma_+}
\big(\langle \zeta, \zeta \rangle -
\langle U_1^{\sigma} \cdots U_N^{\sigma} \ \zeta, U_1^{\sigma} \cdots U_N^{\sigma} \ \zeta \rangle \big) \Big]
\end{equation*}
\begin{equation}\label{T*-N-W}
              \times \   W(U_1^{\sigma} \cdots U_N^{\sigma} \ \zeta) \ ,
\end{equation}
where we denote $U_n^{\sigma} := U_n^{\sigma}(\tau) $.
\end{cor}

To illustrate the above statements by an example, we consider the evolution of the
initial state given by product of the Gibbs states:
\begin{equation}\label{productGibbs}
\rho = \rho_0\otimes\bigotimes_{k=1}^N  \rho_k \ , \
\rho_0= e^{-\beta_0 a^*a}/Z(\beta_0)\  , \ \rho_j = e^{-\beta a^*a}/Z(\beta)  \
( j = 1,2 \ldots, N) ,
\end{equation}
which has the characteristic function (see \cite{TZ}):
\begin{equation}\label{W-L41}
\omega_{\rho}(W(\zeta)) =   \langle \rho \, | \, W(\zeta) \rangle\\
   = \exp\Big[-\frac{|\zeta_0|^2}{4}
         \Big(\frac{1+e^{-\beta_0}}{1-e^{-\beta_0}}
       - \frac{1+e^{-\beta}}{1-e^{-\beta}}\Big)
        - \frac{\langle \zeta,\zeta \rangle}{4}
        \frac{1+e^{-\beta}}{1-e^{-\beta}}\Big] \, .
\end{equation}

From Corollary \ref{cor63}, we obtain the following proposition about time evolution of
the Gibbs state for the open system $(\mathcal{S} + \mathcal{R}) + \mathcal{C}$.
\begin{prop}\label{cor64}
Let $\rho$ be initial density matrix (\ref{productGibbs}).
Then
\[
  \omega_{T_{N\tau,0}^{\sigma}\rho}(W(\zeta))
   = \langle \rho \, | \, T_{N\tau,0}^{\sigma \, *}(W(\zeta)) \rangle
 = \exp\Big[-\frac{1}{4}\langle \zeta, X^{\sigma}(N\tau) \zeta\rangle\Big],
\]
where $ X^{\sigma}(N\tau) $ is the $(N+1)\times(N+1)$ matrix given by
\begin{eqnarray}\label{X-sigma}
 X^{\sigma}(N\tau) &=& U_N^{\sigma \; *} \cdots U_1^{\sigma \; *}
\Big[\Big( - \frac{\sigma_- + \sigma_+}{\sigma_- - \sigma_+} +
\frac{1+e^{-\beta}}{1-e^{-\beta}} \Big) I +
\Big(\frac{1+e^{-\beta_0}}{1-e^{-\beta_0}}-\frac{1+e^{-\beta}}{1-e^{-\beta}}\Big)P_0\Big]
\nonumber
\\
      && \times U_1^{\sigma} \cdots U_N^{\sigma}
      \quad + \ \frac{\sigma_- + \sigma_+}{\sigma_- - \sigma_+}I .
\end{eqnarray}
\end{prop}
\vspace{-0.8cm}
\subsection{{Completely positive quasi-free maps and states}}\label{CP-MDS}
\noindent Let $\mathfrak{A}(\mathfrak{S}, \sigma)$ be the (abstract) Weyl CCR-algebra
for a  linear space $\mathfrak{S}$ and a symplectic form $\sigma$ on $\mathfrak{S}$.
Recall that a bounded linear unital map
$\mathfrak{T}: \mathfrak{A}(\mathfrak{S}, \sigma) \rightarrow \mathfrak{A}(\mathfrak{S}, \sigma)$ is
\textit{quasi-free} if there exists a linear map $U: \mathfrak{S} \rightarrow \mathfrak{S}$,
and a map $\Gamma: \mathfrak{S} \rightarrow \mathbb{C}$  such that
\begin{equation}\label{Q-F-CP}
              \mathfrak{T}(W(\zeta)) = \Gamma(\zeta) W(U \zeta) \quad \
               {\rm{ hold \ for \ all}} \quad \zeta \in\mathfrak{S} \, .
\end{equation}

We also recall that for two $C^*$-algebras $\mathfrak{A}$ and $\mathfrak{B}$, a map $\mathfrak{T}: \mathfrak{A} \rightarrow
\mathfrak{B}$ is \textit{completely positive} (CP) if
\begin{equation}\label{CP}
\sum_{k,k' =1}^{K} y^{*}_k \, \mathfrak{T}(x^{*}_k x_{k'}) \, y_{k'} \geqslant 0
\end{equation}
holds for all $\{x_k\}_{k=1}^{K}\subset \mathfrak{A}$ and
$\{y_k\}_{k=1}^{K}\subset \mathfrak{B}$ for  any  $K \geqslant 1$.
See \cite{Pe} Ch.8, \cite{AF}, \cite{AJP3}.
Using (\ref{Q-F-CP}), one can define the map $\mathfrak{T}$ for a given $U$ and $\Gamma$
on the algebraic span of Weyl operators which is dense in $\mathfrak{A}(\mathfrak{S}, \sigma)$.
For the problem of extension of $\mathfrak{T}$ to a CP map on $\mathfrak{A}(\mathfrak{S}, \sigma)$,
we refer the following result in \cite{DVV}.
\begin{prop}\label{DVV}
For a given linear map $U: \mathfrak{S} \rightarrow \mathfrak{S}$, let $\sigma_U$ be
another symplectic form defined by
\begin{equation}\label{sympl-form-U}
\sigma_{U} (\alpha, \beta) = \sigma (\alpha, \beta) - \sigma (U \alpha, U \beta)  \qquad
 \mbox{for} \quad \alpha, \beta \in \mathfrak{S} \, .
\end{equation}
Then the necessary and sufficient condition of that the map (\ref{Q-F-CP}) can be
extended to a completely positive map on $\mathfrak{A}(\mathfrak{S}, \sigma)$ is the
existence of a state $\omega$ on $\mathfrak{A}(\mathfrak{S}, \sigma_U)$ such that
$\Gamma (\zeta)= \omega(W_{U}(\zeta))$ for Weyl operators $W_{U} \in \mathfrak{A}(\mathfrak{S}, \sigma_U)$.
\end{prop}
\begin{thm}\label{MDS-CP} As a map on $\mathscr{A}(\mathscr{H})$,
the dual MDS $\{{T}_{n,t}^{\sigma \; *}\}_{t\geq 0}$
given by the duality (\ref{T-*sigma}) is quasi-free and completely positive for $n= 1,2, \ldots N$.
\end{thm}
{\sl Proof:} It is obvious from (\ref{T*-n}) and its contraction property that the dual MDS
maps $\mathscr{A}(\mathscr{H})$ into itself and that it is quasi-free.

For a fixed $n$ and $t$, we put $U=  U_n^{\sigma}(t)$.
By setting $\mathfrak{S} = \C^{N+1}$ and
$\sigma( \, \cdot \, , \, \cdot \, ) = {\rm Im} \, \langle  \, \cdot \, , \, \cdot \,  \rangle$,
then $\mathfrak{A}(\mathfrak{S}, \sigma) = \mathscr{A}(\mathscr{H})$ holds and
the action of $T_{n, t}^{\sigma \; *}$ has the form (\ref{Q-F-CP}).
Since $U$ is a contraction (Remark \ref{E-param}), there is a linear map
$C: \C^{N+1} \rightarrow \C^{N+1}$ such that
\begin{equation}\label{C}
\langle C \alpha, C \beta \rangle = \langle \alpha, \beta \rangle - \langle U \alpha, U \beta \rangle \ , \ \ \
\alpha,\beta \in \mathbb{C}^{N+1} \ .
\end{equation}
Then we can consider the CCR-algebra $\mathfrak{A}(\C^{N+1}, \sigma_U)$ with symplectic form
\begin{equation}\label{sigma-U-t}
       \sigma_{U} (\alpha, \beta) = \sigma(\alpha, \beta) - \sigma( U\alpha, U\beta)
           = {\rm{Im}}\langle C \alpha, C \beta \rangle \, ,
\end{equation}
as the $C^*$-subalgebra of $\mathcal{L}(\mathscr{H})$ generated by the Weyl system
$\{W(C \zeta) \mid \zeta \in \mathbb{C}^{N+1}\}$.
Note that $W_{U}(\zeta) = W(C \zeta)$ satisfies CCR with symplectic form $\sigma_U$,
where $W(\zeta)$ is given by (\ref{Wz-bis}).

Let $\rho$ be the product density matrix (\ref{productGibbs}) with
$\beta_0 = \beta = \log\sigma_-/\sigma_+ > 0 $ (c.f. \textbf{H.2}).
Then we have the corresponding normal state $\omega_{\rho}$ on $\mathcal{L}(\mathscr{H})$.
Let $\omega$ be the restriction of $\omega_{\rho}$ to $\mathfrak{A}(\C^{N+1}, \sigma_U)$.
From (\ref{W-L41}), one gets
\[
   \omega(W_U(\zeta)) = \langle \rho \, | \, W(C\zeta)\rangle =
        = \exp \Big[- \frac{\| C \zeta \|^{2})}{4}\,
        \frac{1\sigma_- + \sigma_+}{\sigma_- - \sigma_+} \Big]\, .
\]
Comparing with $\Gamma_{n,t}^{\sigma}(\zeta)$ in (\ref{Omega}) and (\ref{C}),
we see that the $\omega$ plays the role in Proposition \ref{DVV}.
Then, there exists a CP map on $\mathscr{A}(\mathscr{H})$ whose action on the Weyl operators
coincides with that of $T_{n, t}^{\sigma \; *}$.
From the continuity, the coincidence of these maps on $\mathscr{A}(\mathscr{H})$ follows.
Thus, the complete positivity of $T_{n, t}^{\sigma \; *}$ has been proved.
\hfill $\square$

\medskip

Since a composition of quasi-free CP maps is clearly quasi-free and CP,
(\ref{T*-t}) imply
\begin{cor}\label{CP-N}
The dual evolutions $T_{t, 0}^{\sigma \ *}$ is the completely positive quasi-free map on
the Weyl CCR-algebra $\mathscr{A}(\mathscr{H})$ for $t \in [ 0, N\tau)$.
\end{cor}
\begin{cor}\label{CP-L}
The dual evolutions $T_{t, 0}^{\sigma \ *}$ and $T_{n, t}^{\sigma \, *} \ (n = 1, \cdots, N)$
are the completely positive maps on
 $\mathcal{L}(\mathscr{H})$ for $t \in [ 0, N\tau)$.
\end{cor}
{\sl Proof : } For $\mathfrak{T} = T_{n, t}^{\sigma ; *}$,
arbitrarily fixed unit vector $\varphi \in \mathscr{H}$ and $K \in \N$, put
\begin{equation}\label{Phi}
  \Phi( \{A_k\}_{k=1}^{K}, \{B_k\}_{k=1}^{K}) =  (\varphi,
       \sum_{k,k' =1}^{K} B^{*}_k \, \mathfrak{T}(A^{*}_k A_{k'}) \, B_{k'} \varphi) \,,
\end{equation}
where $\{A_k\}_{k=1}^{K}, \{B_k\}_{k=1}^{K} \subset \mathcal{L}(\mathscr{H})$.
Since the CCR algebra $\mathscr{A}(\mathscr{H})$ is a dense subset of $ \mathcal{L}(\mathscr{H}) $
in the strong operator topology,
we may take $\{A_{k,j}\}_{j \in \N}, \{B_{k,j}\}_{j \in \N} \subset \mathscr{A}(\mathscr{H}) $
such that
\[
       {\rm s}\! -\!\! \lim_{j\to \infty} A_{k,j} = A_k  \qquad \mbox{and} \qquad
          {\rm s}\!-\!\!\lim_{j\to \infty}B_{k,j} = B_k
\]
for every $k = 1, \cdots, K$.
Recalling that $\mathfrak{T} = T_{n, t}^{\sigma \, *}$ is CP on $\mathscr{A}(\mathscr{H})$,
we have
\[
    0 \leqslant \Phi( \{A_{k, j}\}_{k=1}^{K}, \{B_{k, j}\}_{k=1}^{K}) =
         \sum_{k, k' =1}^K  \langle T_{n, t}^{\sigma}(B_{k, j}P_{\varphi}
          B_{k', j}^*), A_{k, j}^*A_{k', j} \rangle \, ,
\]
where $P_{\varphi}$ is the projection operator on $\mathscr{H}$ onto its one dimensional
subspace spanned by $\varphi$.
Note that
\[
                  T_{n, t}^{\sigma}(B_{k, j}P_{\varphi} B_{k', j}^* )
             \rightarrow T_{n, t}^{\sigma}(B_kP_{\varphi} B_{k'}^* )  \qquad
         \mbox{in} \quad \mathfrak{C}_1(\mathscr{H})
\]
as $j \to \infty$ since $ B_{k, j}P_{\varphi} B_{k', j}^* \rightarrow
B_kP_{\varphi} B_{k'}^* $ in $ \mathfrak{C}_1(\mathscr{H}) $
and $T_{n, t}^{\sigma}$ is bounded on $ \mathfrak{C}_1(\mathscr{H}) $.
Note also that $A_{k, j}^*A_{k', j}$ converges to $A_k^*A_{k'}$ weakly.
By the principle of uniform boundedness, $\{ A_{k, j}^*A_{k', j}\}_{ j \in \N}$
is a bounded set.
Together with weak continuity of normal states $\langle \rho | \, \cdot \, \rangle:
\mathcal{L}(\mathscr{H}) \rightarrow \C$, this yields that
\[
      \Phi( \{A_k\}_{k=1}^{K}, \{B_k\}_{k=1}^{K}) = \lim_{j\to\infty}
         \Phi( \{A_{k, j}\}_{k=1}^{K}, \{B_{k, j}\}_{k=1}^{K}) \geqslant 0 \, .
\]
Thereby we have proved the  complete positivity of $\mathfrak{T} = T_{n, t}^{\sigma \, *}$
on $\mathcal{L}(\mathscr{H})$.
Proofs for the other maps are almost verbatim.
\hfill$\square$

\medskip

As we have seen, $\mathfrak{T}\big( \mathscr{A}(\mathscr{H}) \big)
\subset \mathscr{A}(\mathscr{H})$
for $\mathfrak{T} = T_{t, 0}^{\sigma \, *}$, $T_{n,t}^{\sigma \, *}$ holds.
Moreover, $\mathfrak{T}$ is  positive unital map.
Therefore for any state $\omega$ on $\mathscr{A}(\mathscr{H})$, $\omega \circ \mathfrak{T}$
is also a state on $\mathscr{A}(\mathscr{H})$.

Recall that a state $\omega$ on $\mathscr{A}(\mathscr{H})$ is said to be quasi-free
if there exist a linear form $L$ and a non-negative sesquilinear form $q$ on $\C^{N+1}$ such that
\[
             \omega(W(\zeta)) = \exp [ L(\zeta) - q(\zeta, \zeta)]
\]
holds for every $\zeta \in \C^{N+1}$\cite{Ve}.
By (\ref{T*-n}) and (\ref{T*-N-W}), it is obvious that
$\omega\circ \mathfrak{T}$ is quasi-free if $\omega$ is.
Let us summarize them in the following proposition.
\begin{prop}\label{QF}
The operators $T_{t, 0}^{\sigma \, *}$, $T_{n, t}^{\sigma \, *} \ (n = 1, \cdots, N)$
map the set of quasi-free state on $\mathscr{A}(\mathscr{H})$
into itself.
\end{prop}

\smallskip

\noindent  \textbf{Acknowledgements }

\noindent
H.T. thanks  JSPS for the financial support under the Grant-in-Aid for Scientific Research (C) 24540168.
He is also grateful to Aix-Marseille and Toulon Universities for their hospitality.
V.A.Z. acknowledges the Institute of Science and Engineering and Graduate School of the Natural Science and Technology
of Kanazawa University for financial support and hospitality allowed to finish this paper.
\section{Appendix}\label{App}
{\sl Proof:} (of Lemma \ref{m-accr})
The operator $K_0$ with its domain $\mathcal{D}_0$ (\ref{domH}) is closed
with discrete spectrum $\mathfrak{S}(K_0) \subset \C_+ = \{ \, z \in \C | \, {\rm{Re}}\,\, z > 0 \}$.
In fact, for $m \in \Z^{N+1}$,  $\big(2^{-1}(\sigma_+ + \sigma_-)+iE\big) m_0 +
2^{-1}\sigma_+ +i\epsilon\sum_{j=1}^N  m_j$ is its eigenvalue and $\psi_m$ in (\ref{CONS})
is the corresponding eigenvectors.
It is enough to consider the case $n =1 $ only and to prove the following three claims \cite{Za}:

\noindent(i) the operator $K_1$ is closed ;

\noindent(ii) the \textit{numerical range} of  $K_1$ is contained in $\overline{\C}_{+}$ ;

\noindent(iii) there exists $z\in\C$ such that  ${\rm{Re}}\, z > 0$ and
$z$ belongs to the \textit{resolvent set}
$\rho(- K_1)$ of the operator $- K_1$.

For (i), we show that there exist constants $c\in (0,1)$ and $C >0$ such that
\begin{equation}\label{cond-close}
 \| \eta(b_0^*b_1 +b_1^*b_0) \varphi \| \leqslant c \| K_0\varphi\| + C \|\varphi\|
\end{equation}
for every $\varphi \in \mathcal{D}_0$.
It is obvious from CCR that
\[
    \|b_0^*b_1 \varphi\|^2 = \big( b_0b_0^* \varphi,
     b_1^*b_1  \varphi \big) \leqslant \big( (b_0^*b_0 + \mathbb{1}) \varphi,
     (b_1^*b_1 + \cdots + b_N^*b_N + \mathbb{1}) \varphi \big),
\]
\[
    \|b_1^*b_0 \varphi\|^2 = \big( b_0^*b_0 \varphi,
      b_1b_1^*  \varphi \big) \leqslant \big( (b_0^*b_0 + \mathbb{1}) \varphi,
     (b_1^*b_1 + \cdots + b_N^*b_N + \mathbb{1}) \varphi \big)
\]
hold for any $\varphi \in \mathcal{D}_0$. Note that
\[
      \Big\|\Big( E -i\frac{\sigma_++\sigma_-}{2}\Big)(b_0^*b_0+ \mathbb{1}) \varphi +
    \epsilon (b_1^*b_1 + \cdots + b_N^*b_N + \mathbb{1}) \varphi \Big\|^2
\]
\[
   - 2 \Big(E +\sqrt{E^2 + {(\sigma_++\sigma_-)^2}/{4}} \Big) \, \epsilon \,
      \big( (b_0^*b_0+ \mathbb{1})  \varphi,  \,
    (b_1^*b_1 + \cdots + b_N^*b_N + \mathbb{1}) \varphi \big)
\]
\[
        = \Big\| \sqrt{E^2 +{(\sigma_++\sigma_-)^2}/{4}} \, (b_0^*b_0+ \mathbb{1})
\varphi  - \epsilon (b_1^*b_1 + \cdots + b_N^*b_N + \mathbb{1}) \varphi \Big\|^2
         \geqslant 0
\]
and
\[
        -iK_0 \varphi + \Big(E+ \epsilon -i\frac{\sigma_-}{2}\Big) \varphi = \Big( E -i\frac{\sigma_++\sigma_-}{2}\Big)
           (b_0^*b_0+ \mathbb{1}) \varphi +
    \epsilon (b_1^*b_1 + \cdots + b_N^*b_N + \mathbb{1}) \varphi \, .
\]
Then, we have
\[
        \|\eta(b_0^*b_1+ b_1^*b_0) \varphi \| \leqslant 2\eta
       \Big(E +\sqrt{E^2 +{(\sigma_++\sigma_-)^2}/{4}}
      \Big)^{-1/2} \epsilon^{-1/2}
\]
\[
      \times \sqrt{\Big( E +\sqrt{E^2 +
     {(\sigma_++\sigma_-)^2}/{4}} \Big)\epsilon   \Big( (b_0^*b_0+ \mathbb{1})  \varphi,
     (b_1^*b_1 + \cdots + b_N^*b_N + \mathbb{1}) \varphi \Big)}
\]
\[
        \leqslant \sqrt 2\eta
       \Big(E +\sqrt{E^2 +{(\sigma_++\sigma_-)^2}/{4}}
      \Big)^{-1/2}\epsilon^{-1/2}\Big\|  -iK_0 \varphi + \Big(E+ \epsilon -i
           \frac{\sigma_-}{2}\Big) \varphi \Big\|
\]
\[
        \leqslant c\|K_0\varphi \| + C\| \varphi \|,
\]
where
\[
       c = \sqrt 2\eta \Big(E +\sqrt{E^2
     + {(\sigma_++\sigma_-)^2}/{4}}\Big)^{-1/2} \epsilon^{-1/2}  <1
\]
because of the conditions $\eta^2 \leqslant E\epsilon $ \textbf{(H.1)} and
$0 < \sigma_+ + \sigma_- $ \textbf{(H.2)}.

To show (ii), let $\varphi \in \mathcal{D}_0$, $\|\varphi\| =1$. Then one gets
\[
        ( \varphi, K_1 \varphi) = \frac{\sigma_+}{2} ( \varphi, b_0b_0^*\varphi)\, +
        \frac{\sigma_-}{2} ( \varphi, b_0^*b_0 \varphi )\,
\]
\[
      + i \, E ( \varphi, b_0^*b_0\varphi) +i\epsilon  ( \varphi, \sum_{j=1}^Nb_j^*b_j \varphi )
       +i \, \eta ( \varphi, (b_0^*b_1 + b_1^*b_0)\varphi)\,
\]
\[
     = \frac{\sigma_+}{2}  \|b_0^*\varphi\|^2
       + \frac{\sigma_-}{2} \|b_0 \varphi \|^2
         + iE  \|b_0\varphi\|^2 + i\epsilon \sum_{j=1}^N\|b_j \varphi \|^2
      +2i \eta \ {\rm{Re}}\, ( b_0 \varphi, b_1\varphi)  \subset \overline{\C}_+  .
\]

For (iii), we note that by virtue of $\mathfrak{S}(K_0) \subset \C_+$,
$z \in \rho(- K_0)$ and $\|(z\mathbb{1} + K_0)^{-1}\| \leqslant 1/{\rm{Re}}\, z$ hold, if ${\rm{Re}}\, z > 0$.
Moreover, the identity
\[
    \| (z\mathbb{1} +K_0)\varphi \|^2 - \|K_0\varphi \|^2 = |z|^2 \|\varphi \|^2
    + 2({\rm{Re}}\, z)\Big( \varphi, \Big( \frac{\sigma_+}{2}  b_0b_0^* +
        \frac{\sigma_-}{2}  b_0^*b_0\Big) \varphi \Big)
\]
\[
       + 2({\rm{Im}}\, z)\Big( \varphi, \Big( Eb_0^*b_0 +
      \epsilon \sum_{j=1}^Nb_j^*b_j \Big) \varphi \Big)
\]
yields $\|K_0(z \mathbb{1} +K_0)^{-1}\| \leqslant 1$, if ${\rm{Re}}\, z \geqslant 0$ and ${\rm{Im}}\, z \geqslant 0$ hold.
Hence, by (\ref{cond-close}), we obtain
\begin{align}\label{estim}
&\| \eta(b_0^*b_1 +b_1^*b_0) (z\mathbb{1} +K_0)^{-1}\|  \nonumber \\
&  \leqslant  c \| K_0(z\mathbb{1} +K_0)^{-1}\|  + C \|(z\mathbb{1} +K_0)^{-1}\| \leqslant
c + \frac{C}{{\rm{Re}}\, z}
\end{align}
for such value of $z$.
Then, if ${\rm{Re}}\, z$ is large enough, the right hand side of (\ref{estim}) is less than one.
For this $z$, thanks to the resolvent identity for $K_1$ and $K_0$, we have the boundedness of
\[
      (z \mathbb{1} +K_1)^{-1} =  (z\mathbb{1} +K_0)^{-1}
     \big(\mathbb{1} + i\eta (b_0^*b_1 +b_1^*b_0)
      (z\mathbb{1} +K_0)^{-1} \big)^{-1} \ ,
\]
which proves the assertion (iii) and the lemma.  \hfill $\square$

\end{document}